# Three-Level Parallel $J$-Jacobi Algorithms for Hermitian Matrices[*]

Sanja Singer[†], Saša Singer[‡], Vedran Novaković[§],
Davor Davidović[¶], Krešimir Bokulić[‖] Aleksandar Ušćumlić[**]

November 2, 2018

**Abstract**

The paper describes several efficient parallel implementations of the one-sided hyperbolic Jacobi-type algorithm for computing eigenvalues and eigenvectors of Hermitian matrices. By appropriate blocking of the algorithms an almost ideal load balancing between all available processors/cores is obtained. A similar blocking technique can be used to exploit local cache memory of each processor to further speed up the process. Due to diversity of modern computer architectures, each of the algorithms described here may be the method of choice for a particular hardware and a given matrix size. All proposed block algorithms compute the eigenvalues with relative accuracy similar to the original non-blocked Jacobi algorithm.

## 1 Introduction

The Jacobi method is one of the simplest algorithms for diagonalization of symmetric matrices. For many years, during the era of predominantly sequential computing, the Jacobi-type algorithms were almost forgotten in practice, as they are too slow for serial computing, when compared with some faster diagonalization algorithms proposed since.

But, in 1992, Demmel and Veselić's historic paper [11] shed a new light onto the method. They proved that the Jacobi algorithm for positive definite matrices is accurate in the relative sense, which is the natural error model for standard floating–point computations.

Let $A$ be a symmetric positive definite matrix of order $n$ with elements $a_{ij}$. If we introduce a perturbation $\delta a_{ij}$ in each element of $A$, such that $\delta a_{ij}$ is small in the relative sense with respect to $a_{ij}$, i.e.,

$$|\delta a_{ij}| \leq \varepsilon |a_{ij}|,$$

[*]This work was supported by grant 037–1193086–2771 by Ministry of Science, Education and Sports, Croatia. A part of this paper was done while the last three authors were students at the Department of Mathematics, University of Zagreb, Croatia.

[†]Faculty of Mechanical Engineering and Naval Architecture, University of Zagreb, Ivana Lučića 5, 10000, Croatia, e-mail: ssinger@fsb.hr

[‡]Department of Mathematics, University of Zagreb, Bijenička cesta 30, 10000 Zagreb, Croatia, e-mail: singer@math.hr

[§]Faculty of Mechanical Engineering and Naval Architecture, University of Zagreb, Ivana Lučića 5, 10000 Zagreb, Croatia, e-mail: venovako@fsb.hr

[¶]Ruđer Bošković Institute, Bijenička 54, 10000 Zagreb, Croatia, e-mail: davor.davidovic@irb.hr

[‖]Multicom d.o.o., V. Preloga 11, 10000 Zagreb, Croatia, e-mail: kresimir.bokulic@multicom-is.hr

[**]Ph.D. student at the Department of Mathematics, University of Zagreb, Bijenička cesta 30, 10000 Zagreb, Croatia, e-mail: sasa.uscumlic@gmail.com

then the eigenvalues $\lambda_i$ of $A$, and $\lambda'_i = \lambda_i + \delta\lambda_i$ of $A + \delta A$, satisfy

$$\frac{|\delta\lambda_i|}{\lambda_i} \leq \varepsilon n \kappa(A_s), \quad i = 1, \ldots, n.$$

Here, $A_s = (\text{diag}(A))^{-1/2} A (\text{diag}(A))^{-1/2}$ is a special diagonal scaling of $A$, where $\text{diag}(A)$ denotes a diagonal matrix containing the diagonal elements of $A$, and $\kappa(A_s) = \|A_s\|_2 \|A_s^{-1}\|_2$ is the scaled condition of $A$. It can also be shown that $\kappa(A_s)$ nearly minimizes the condition number over all possible diagonal scalings of $A$ (see [33] for details).

Almost at the same time, Veselić [34] has developed a Jacobi-type algorithm (sometimes called indefinite, hyperbolic or $J$-Jacobi algorithm) for diagonalization of definite matrix pairs $(A, J)$, where $A$ is symmetric positive definite, and $J$ is a diagonal matrix of signs 1 and $-1$. In contrast to the ordinary Jacobi method, which uses unitary congruences, the $J$-Jacobi method uses $J$-unitary congruences for diagonalization. This algorithm can also be used for diagonalization of indefinite symmetric matrices $H$ (see Section 2). As a Ph.D. student of Veselić, Slapničar [30, 32] has shown that the algorithm is accurate in the relative sense, as well. The same is true for various implementations of the $J$-Jacobi algorithm for complex Hermitian matrices [27].

Recently, Dopico, Koev and Molera [12] have shown that the ordinary one-sided Jacobi algorithm can also be relatively accurate, if a symmetric, possibly indefinite, matrix $H$ can be factorized as $H = XDX^T$, where $D$ is diagonal, and the factor $X$ is well-conditioned. But, according to their own tests, the ordinary one-sided algorithm can have significantly more sweeps than the one-sided $J$-Jacobi algorithm. Moreover, when the eigenvectors of $H$ are required, the $J$-Jacobi algorithm is much faster, even if both algorithms have the same number of sweeps. This is due to the fact that, after completion of the $J$-Jacobi process, the eigenvectors are obtained by simple scaling of the final matrix, while in the ordinary algorithm the eigenvectors are obtained by accumulation of transformations. Therefore, the $J$-Jacobi algorithm seems to be the method of choice for accurate computation of eigenvalues and eigenvectors of symmetric indefinite matrices.

It is well-known that Jacobi-type algorithms are almost ideally parallelizable. Many authors have published, mostly two-sided, parallel implementations of the ordinary algorithm (see, for example, [1, 3, 5, 15, 22, 24, 35, 37]). Convergence results for different parallel orderings in the ordinary Jacobi-type algorithms are given in [13, 20, 21, 26].

Besides the inherent element-wise parallelism, the blocking paradigm can also be exploited to enhance the performance, especially in one-sided versions of the algorithms. In practice, the idea of blocking can be applied with two different goals in mind—to achieve data independency for parallelization, or to reuse large chunks of data. These two opposite goals for blocking can be used at different levels of the algorithm to increase its efficiency.

Blocking to reuse data can be used in sequential Jacobi-type algorithms to achieve significant speedups in presence of a two-level memory hierarchy on a single processor (see [17, 18]). Quite generally, there are two different approaches to blocking of these algorithms: the *block-oriented* approach (see [17]), and the *full block* approach (see [18]). For some well-known classes of pivoting strategies, including the block column-cyclic strategy, both approaches lead to globally convergent algorithms [16, 17, 18]. These blocked versions of the algorithms retain essentially the same accuracy as the non-blocked ones [17].

In this paper, we present several three-level and two-level parallel implementations of the one-sided $J$-Jacobi algorithm. In three-level implementations, the *outer* level of blocking is used for parallelization of the algorithm, while the *inner* level of blocking targets the local memory hierarchy. Finally, the bottom level of actual computation (in terms of elements of matrices) efficiently uses the local cache memory of each processor/core. The inner level of blocking may be omitted if the working blocks of matrices are small enough to reside in cache throughout the whole computation phase. These simplified two-level parallel implementations may be significantly more efficient for smaller matrices.

For larger matrices, the inner level of blocking in each algorithm brings additional speedup over the two-level parallel algorithm. This speedup depends on the algorithm, the matrix size, the number of processes used, and the interprocess communication involved. For example, for matrices of order $n = 16000$ on 16 processor cores, the speedup for all used three-level versions of the algorithm is over 40%, and may reach almost 70%. The crossover point in efficiency between two-level and three-level implementations depends heavily on a particular hardware architecture, and should be determined by numerical testing.

The global convergence of parallel algorithms can also be proved, if we use the block modulus pivot strategy (with slight modifications for inclusion of the diagonal blocks) for the outer level of blocking. This follows from the fact that the block modulus pivot strategy is weakly equivalent to the block antidiagonal ordering, and this strategy is weakly equivalent to the block row-cyclic strategy (see [26]). Finally, block-cyclic strategy (with cyclic by rows or columns inner strategy) is convergent (see [17, 18]).

The paper is organized as follows. First, we give a detailed description of the one-sided $J$-Jacobi algorithm. Section 3 deals with sequential blocked versions of the algorithm. In Sections 4 and 5 we describe the two-level and three-level parallel implementations of blocked algorithms, respectively. The final section contains the results of numerical testing.

## 2 Two-sided and one-sided Jacobi algorithms

In this section we shall establish a connection between two-sided and one-sided versions of the Jacobi algorithm. We begin by describing the ordinary two-sided Jacobi algorithm, as it is easier to understand. Then we derive the one-sided hyperbolic algorithm and describe some of its advantages for practical computation.

### 2.1 The two-sided Jacobi algorithm

The ordinary two-sided Jacobi algorithm for diagonalization of a given Hermitian matrix $H$ performs a sequence of similarity transformations. Each transformation annihilates a single off-diagonal element in the working matrix (see [23, 25]). The similarity used to annihilate the element in position $(r, s)$ consists of a trigonometric rotation $U_T(r, s)$ in the $(r, s)$ plane. In the real case, it is equal to the identity matrix, except at the intersections of the $r$-th and $s$-th rows and columns, where, in Matlab notation,

$$U_T(r,s)([r,s],[r,s]) = U_T := \begin{bmatrix} \cos\varphi & \sin\varphi \\ -\sin\varphi & \cos\varphi \end{bmatrix}. \tag{2.1}$$

In the complex case, we can use various forms of trigonometric plane rotations $U_T$. For example,

$$U_T = \begin{bmatrix} \cos\varphi & e^{i\alpha}\sin\varphi \\ -e^{-i\alpha}\sin\varphi & \cos\varphi \end{bmatrix}, \tag{2.2}$$

or

$$U_T = \begin{bmatrix} e^{i\alpha}\cos\varphi & e^{i\alpha}\sin\varphi \\ -\sin\varphi & \cos\varphi \end{bmatrix}, \quad \text{or} \quad U_T = \begin{bmatrix} \cos\varphi & \sin\varphi \\ -e^{-i\alpha}\sin\varphi & e^{-i\alpha}\cos\varphi \end{bmatrix}. \tag{2.3}$$

Since $U_T(r,s)$ is a unitary matrix, each similarity transformation is also a unitary congruence. The whole sequence of transformations is usually divided into sweeps. In each sweep, all $n(n-1)/2$ off-diagonal elements in, say, the upper triangle of the working matrix are annihilated exactly once, according to some prescribed traversing or annihiliation order, also known as the pivot strategy. The following strategies are in widespread use:
 1. row- and column-cyclic strategies, with the advantage of sequential data access (thus, well suited for the inner level of blocking),
 2. modulus and round-robin strategies, frequently used for parallelization.

Quite generally, the two-sided Jacobi process on an initial matrix $H^{(0)} = H$ is an iterative process of the following form:

$$H^{(S)} = [U^{(S)}]^* \cdots [U^{(2)}]^*[U^{(1)}]^* H U^{(1)} U^{(2)} \cdots U^{(S)}, \qquad (2.4)$$

where $U^{(i)}$ denotes a product of all $n(n-1)/2$ trigonometric rotations applied in sweep $i$, for $i = 1, \ldots, S$. The process terminates after, say, $S$ sweeps, when all off-diagonal elements of $H^{(S)}$ are relatively small compared to $\mathrm{diag}(H^{(S)})$. Then, the diagonal elements of $H^{(S)}$ are the computed eigenvalues, and $U := U^{(1)} U^{(2)} \cdots U^{(S)}$ is the computed eigenvector matrix.

## 2.2 Factorization for the one-sided algorithm

The iterative process (2.4) can also be organized as a one-sided algorithm. In the first phase, the method computes the Hermitian indefinite factorization of the initial matrix $H$, by a variant of the Bunch–Parlett factorization (see [2, 6, 7, 8, 9, 10])

$$PHP^T = MDM^*.$$

Here, $P$ is a permutation matrix, $M$ is a lower triangular matrix, and $D = D^*$ is a block-diagonal matrix with diagonal blocks of order 1 or 2. An additional diagonalization of the diagonal blocks of $D$, and an appropriate scaling of the columns of $M$ (see [31] for details) yield

$$PHP^T = GJG^*, \qquad J = \mathrm{diag}(j_{11}, \ldots, j_{nn}). \qquad (2.5)$$

If $H$ is nonsingular, then $G$ is lower block-triangular with diagonal blocks of order 1 or 2, and $J$ contains the inertia of $H$, i.e., $j_{ii} \in \{-1, 1\}$, for $i = 1, \ldots, n$. If $H$ is singular and $\mathrm{rank}(H) = m < n$, then $G$ in (2.5) is a lower block-trapezoidal $n \times m$ matrix of full column rank. In this case, $J$ is of order $m$ and contains the signs of nonzero eigenvalues of $H$. Note that $G^*G$ is always positive definite, which will be sufficient for the one-sided $J$-Jacobi algorithm (see below).

For a nonsingular matrix $H$, the *ordinary* one-sided Jacobi algorithm applies, from the right on $G^*$, the sequence of rotations that diagonalizes $H$ (or $PHP^T$) in (2.4). This algorithm can be accurate in the relative sense [12], but is significantly less efficient than its hyperbolic counterpart.

The *hyperbolic* or $J$-Jacobi algorithm for the eigenproblem $Hx = \lambda x$ is obtained by using the factored form (2.5) of $H$,

$$GJG^*x = \lambda x.$$

Pre-multiplication by $G^*$ and $J = J^{-1}$ give

$$G^*Gz = \lambda J z, \quad z = JG^*x. \qquad (2.6)$$

If $G$ is of full column rank, then the eigenproblem for $H$ (except for zero eigenvalues) is equivalent to the eigenproblem for the pair $(G^*G, J)$ from (2.6), with positive definite matrix $A := G^*G$. Thus, we can apply the Jacobi eigenreduction algorithm for positive definite pairs by Veselić [34]. If $H$ is singular, the eigenvectors corresponding to the zero eigenvalue can be computed by a post-processing step, described below.

The matrix $G$ in (2.6) does not have to be computed by the Bunch-Parlett factorization. In fact, it could be any full column rank factor of $H$, such that (2.5) holds. This is useful in some applications where $H$ is given implicitly by (2.5), and $G$ is readily available.

## 2.3 The one-sided $J$-Jacobi algorithm

In the second phase, the eigenvalues of the pair $(A, J)$ are computed by simultaneous diagonalization of $A$ and $J$. This is done by a sequence of $J$-unitary congruences, i.e., matrices $V$ such that $V^* J V = J$. These transformations preserve the structure of $J$, and only $A$ is actually diagonalized in the process. Since $A = G^* G$, this is equivalent to the orthogonalization of columns of $G$, and will be used to compute the eigenvectors of $H$.

Simple $J$-unitary matrices are either trigonometric, or hyperbolic plane rotations. In the real case, a hyperbolic plane rotation $U_H(r, s)$ is equal to identity matrix, except at the intersections of the $r$-th and $s$-th rows and columns, where

$$U_H(r, s)([r, s], [r, s]) = U_H := \begin{bmatrix} \cosh \varphi & \sinh \varphi \\ \sinh \varphi & \cosh \varphi \end{bmatrix}. \tag{2.7}$$

In the complex case, there are several variants of hyperbolic rotations $U_H$, similar to the trigonometric case (2.2)–(2.3). The trigonometric rotation $U_T(r, s)$ is $J$-unitary if $j_{rr} = j_{ss}$, while the hyperbolic rotation $U_H(r, s)$ is $J$-unitary if $j_{rr} = -j_{ss}$. In Jacobi algorithms, the corresponding $J$-unitary congruence is used to annihilate the $(r, s)$ element in the working matrix $A$.

The hyperbolic or $J$-Jacobi algorithm systematically annihilates the elements of the working matrix $A$, until it is diagonalized to a desired precision. As before, we have two-sided and one-sided versions of the algorithm.

Let $G^{(0)} = G$ be the initial factor of $H$ from (2.5), and let $A^{(0)} = A = G^* G$ be the corresponding initial positive definite matrix in (2.6).

To simplify the notation, regardless of the chosen pivoting strategy, let $W^{(k)}$ be the plane rotation used in the $k$-th step of the process. If $(r, s)$ denotes the pair of indices determined by the pivoting strategy in step $k$, then $W^{(k)}$ is either $U_T(r, s)$, or $U_H(r, s)$.

In *two-sided* $J$-Jacobi algorithms, $J$-unitary matrices $W^{(k)}$ are applied on *both* sides of $A$, and the transformation in step $k$ can be written as

$$A^{(k)} = [W^{(k)}]^* A^{(k-1)} W^{(k)}, \qquad J^{(k)} = [W^{(k)}]^* J^{(k-1)} W^{(k)} = J. \tag{2.8}$$

On the other hand, in *one-sided* algorithms, they are applied from the *right* on $G$, and the transformation in step $k$ is

$$G^{(k)} = G^{(k-1)} W^{(k)}. \tag{2.9}$$

Since, initially, $A = G^* G$ and $G$ is of full column rank, it follows immediately that in each step $A^{(k)} = [G^{(k)}]^* G^{(k)}$ and $G^{(k)}$ is also of full column rank, so $A^{(k)}$ remains positive definite throughout the process. The annihilation of element $a_{rs}^{(k-1)}$ in (2.8) is equivalent to the mutual orthogonalization of two pivot columns $g_r^{(k-1)}$, $g_s^{(k-1)}$ of $G^{(k-1)}$ in (2.9).

Numerically stable formulae for computing the elements of $W^{(k)}$, and the transformed elements in the new working matrix $A^{(k)}$ or $G^{(k)}$, can be found in [34].

The algorithm stops after, say, $l$ steps, when the working matrix $A^{(l)}$ becomes diagonal to the working precision. Then, in the one-sided algorithm, the columns of the final matrix $\widetilde{G} := G^{(l)}$ are orthogonal to the working precision. Therefore, $\widetilde{G}$ can be written as

$$\widetilde{G} = \widetilde{U} \Sigma, \tag{2.10}$$

where $\widetilde{U}$ has orthonormal columns, and $\Sigma$ is diagonal, containing the norms of the columns of $\widetilde{G}$. Both $\Sigma$ and $\widetilde{U}$ can be computed trivially from the final $\widetilde{G}$, and the columns of $\widetilde{U}$ are orthonormal to the working precision. Let

$$W := W^{(1)} W^{(2)} \cdots W^{(l)}$$

be the product of all used $J$-unitary matrices that diagonalizes the pair $(A, J)$, and let $V = W^{-*}$. Then, $W$ and $V$ are $J$-unitary, as well, since $J$-unitary matrices form a multiplicative group. From (2.9) we see that $\widetilde{G} = GW$, and (2.10) yields

$$\widetilde{G} W^{-1} = \widetilde{U} \Sigma W^{-1} = \widetilde{U} \Sigma V^* = G. \tag{2.11}$$

The last equality is actually the hyperbolic singular value decomposition (HSVD) of $G$ (see [4]), computed by the one-sided $J$-Jacobi algorithm of Veselić [34].

The general form of the HSVD of a rectangular $n \times m$ matrix $G$ with respect to $J$ is described in [36]. When $n \geq m$ and $\operatorname{rank}(GJG^*) = \operatorname{rank}(G) = m$, it has the following form

$$G = U \begin{bmatrix} \Sigma \\ 0 \end{bmatrix} V^*, \qquad \Sigma = \operatorname{diag}(\sigma_1, \sigma_2, \ldots, \sigma_m), \quad \sigma_1 \geq \sigma_2 \geq \cdots \geq \sigma_m > 0,$$

where $U$ is unitary of order $n$, $\Sigma$ is a diagonal matrix of positive hyperbolic singular values, and $V$ is $J$-unitary of order $m$. If $H = GJG^*$, resembling (2.5), then, by using the HSVD of $G$, we obtain

$$H = GJG^* = U \begin{bmatrix} \Sigma \\ 0 \end{bmatrix} V^* JV \begin{bmatrix} \Sigma & 0 \end{bmatrix} U^* = U \begin{bmatrix} J\Sigma^2 & 0 \\ 0 & 0 \end{bmatrix} U^*.$$

The squares of the hyperbolic singular values of $G$ are, up to the signs in $J$, the nonzero eigenvalues of $H$, and $U$ is the corresponding eigenvector matrix.

If $H$ is singular, of rank $m < n$, then $G$ and $\widetilde{U}$ in (2.11) are rectangular $n \times m$ matrices, and $\widetilde{U}$ contains the eigenvectors that correspond to the nonzero eigenvalues of $H$. Once we compute $\widetilde{U}$, the remaining eigenvectors for the zero eigenvalue span the orthogonal complement of the column space of $\widetilde{U}$. If required, they can be computed by the full QR factorization of $\widetilde{U}$ (see [19]).

Note that the one-sided $J$-Jacobi algorithm *directly* computes the eigenvector matrix $U$, and numerical orthogonality of the computed eigenvectors is guaranteed by the termination test for the process (2.9). There are several other advantages of the one-sided algorithm over the two-sided.

(a) Even if the elements of $W^{(k)}$ in (2.8) and (2.9) are not computed very accurately, this will not harm the one-sided process, as long as $W^{(k)}$ is numerically $J$-unitary to the working precision. The one-sided algorithm will then reduce the off-diagonal norm of the *pivot* submatrix

$$A_P^{(k-1)} := \begin{bmatrix} a_{rr}^{(k-1)} & a_{rs}^{(k-1)} \\ \bar{a}_{rs}^{(k-1)} & a_{ss}^{(k-1)} \end{bmatrix} = \left[ g_r^{(k-1)}, g_s^{(k-1)} \right]^* \left[ g_r^{(k-1)}, g_s^{(k-1)} \right],$$

instead of diagonalizing it. On the other hand, in the two-sided algorithm, the new off-diagonal element $a_{rs}^{(k)}$ is explicitly set to zero, so the computed elements of $W^{(k)}$ should be as accurate as possible.

(b) The two-sided algorithm is much slower than the one-sided algorithm. In each step, the one-sided algorithm changes only the two pivot columns of the working matrix $G^{(k-1)}$,

$$\left[ g_r^{(k)}, g_s^{(k)} \right] = \left[ g_r^{(k-1)}, g_s^{(k-1)} \right] W_P^{(k)}, \tag{2.12}$$

where $W_P^{(k)}$ is of order 2, and denotes the nontrivial, pivotal part of $W^{(k)}$, i.e., $W_P^{(k)}$ is either $U_T$ from (2.1)–(2.3), or $U_H$ from (2.7), including the complex forms of $U_H$. This operation—the so-called *update* of columns, is easily vectorized by modern compilers. On the contrary, the two-sided algorithm changes some parts of two rows and two columns in one triangle of $A^{(k-1)}$. Vectorization and other optimizations of such an operation still impose great challenges to present-day compilers.

Henceforward, we shall often omit the superscripts when referring to the parts of $G^{(k-1)}$, indicating that the parts involved refer to the *current* state of the working matrix $G$. For example, the update of columns (2.12) will be denoted simply by

$$[g_r', g_s'] = [g_r, g_s] W_P^{(k)}, \tag{2.13}$$

where primes denote the new (updated) values of columns.

Finally, numerical testing shows that the $J$-Jacobi method becomes faster by approximately 10%, if $J$ is ordered as $J = \mathrm{diag}(I_\nu, -I_{n-\nu})$, where $\nu \in \{0, \ldots, n\}$ is determined by the inertia of $H$. This ordering can be achieved by using the congruence transformation $(G^*G, J) \mapsto (P_1^* G^* G P_1, P_1^* J P_1)$, with a suitable permutation matrix $P_1$. From now on, we assume that $J$ is already partitioned in such a way.

The optimal choice of the pivoting strategy for the non-blocked version of the one-sided algorithm depends on how matrices are stored within a particular programming environment. Since all our test programs are written in Fortran, where matrices are stored in column-major order, we use the *column-cyclic* strategy for non-blocked algorithms.

## 3 Sequential blocked versions of $J$-Jacobi algorithm

On hardware architectures with two or more layers of memory with different bandwidths, many blocked algorithms show significant speedups. The same can be done by blocking in the second, iterative part of one-sided Jacobi algorithms.

Two blocked modifications of the one-sided $J$-Jacobi algorithm have been proposed in [17, 18]: the *block-oriented* algorithm and the *full block* algorithm. The aim of blocking in these algorithms is to speed-up the most expensive operation, and that is the update of columns (2.13) in the working matrix $G$, by turning it from a BLAS 1 into a BLAS 3 operation. Therefore, both algorithms use a one-dimensional partition of $G$ into blocks of full-length columns.

### 3.1 Block-column partition

A blocked version of the algorithm begins by a block-column partition of the factor $G$ into $b$ blocks, and the corresponding partition of the diagonal matrix $J$

$$G = [G_1, \ldots, G_b], \quad J = \mathrm{diag}(J_{11}, \ldots, J_{bb}), \qquad (3.1)$$

where the block-column $G_\ell$ has $n_\ell$ columns, and the diagonal square block $J_{\ell\ell}$ is of order $n_\ell$, for $\ell = 1, \ldots, b$.

For sequential blocking in practice, the number of blocks $b$ is calculated from the requirement that all block-sizes $n_\ell$ should be as close as possible to a given *target* block-size $n_t$, without exceeding it. The optimal choice of $n_t$ for each algorithm should be determined by numerical testing, and will be illustrated at the end of this section. Once we know the value of $n_t$, the number of blocks is given by

$$b = \left\lceil \frac{n}{n_t} \right\rceil = \left\lfloor \frac{n + n_t - 1}{n_t} \right\rfloor. \qquad (3.2)$$

Then we divide $G$ into $b$ blocks with almost the same number of columns. There are two easy ways to determine the actual block-sizes. Just for simplicity, we shall assume that the block-sizes are decreasingly ordered $n_1 \geq \cdots \geq n_b$.

(a) In a *greedy* partition, all the blocks, except possibly the last one, have the same maximal allowed size $n_t$. The last block may have a smaller number of columns, given by

$$n_b = n - (b-1)n_t = (n-1) \bmod n_t + 1.$$

This is also valid when $n_t \geq n$, i.e., when $b = 1$ in (3.2). Note that if $b > 1$, the last block can have a significantly smaller number of columns than all the other blocks, for $n_1 - n_b$ can be as large as $n_t - 1$.

(b) In a *uniform* partition, this cannot happen, since we require that $n_1 - n_b \leq 1$, so the numbers of columns in any two blocks may differ by at most one.

Because of that, the uniform partition is more suitable for blocked parallel implementations of the algorithms, as it ensures a better load balancing between processes. Consequently, to unify the algorithms, it will be used for all types and levels of blocking, from now on.

For a given number of blocks $b$, let $b_r = n \bmod b$, and $n_{\min} = \lfloor n/b \rfloor$. The block-sizes in the uniform partition are then given by

$$n_\ell = \begin{cases} n_{\min} + 1, & \ell = 1, \ldots, b_r, \\ n_{\min}, & \ell = b_r + 1, \ldots, b. \end{cases} \qquad (3.3)$$

This follows easily from

$$\sum_{\ell=1}^{b} n_\ell = n = \left\lfloor \frac{n}{b} \right\rfloor \cdot b + b_r = \left( \left\lfloor \frac{n}{b} \right\rfloor + 1 \right) \cdot b_r + \left\lfloor \frac{n}{b} \right\rfloor \cdot (b - b_r).$$

Here we can have $b_r = 0$, and then all blocks have the same block-size $n_{\min}$. Otherwise, the maximal block-size is $n_{\min} + 1$. In contrast to the greedy partition, the maximal block-size $n_1$ can be strictly smaller than the target block-size $n_t$, even if $b > 1$.

## 3.2 Block-pivot submatrices

All blocked versions of Jacobi algorithms have a global structure very similar to the non-blocked version. The computation proceeds in block-steps that are usually grouped into block-sweeps, controlled by some prescribed block-pivoting strategy. This strategy determines the blocks that are transformed in each block-step.

For sequential blocked algorithms we shall use the *block column-cyclic* pivoting strategy at the block level. The same strategy will also be used in three-level parallel algorithms for the inner level of blocking, inside each process.

In a blocked version of the one-sided $J$-Jacobi algorithm, each block-step performs a certain $J$-unitary transformation of the working matrix $G$. This transformation has the same form as in (2.9),

$$G^{(k)} = G^{(k-1)} W^{(k)}, \qquad (3.4)$$

where $k$ now refers to a block-step, instead of a single step. Various versions differ only in the choice of $W^{(k)}$, to be described below.

To begin with, all matrices $G^{(k)}$ in (3.4) are partitioned in the same way as the initial matrix $G^{(0)} = G$, so the partitioning of $G^{(k)}$ does *not* depend on $k$. Moreover, any column partition of $G$ in (3.1) induces a two-dimensional partition of $A = G^*G$ into $b^2$ blocks

$$A_{ij} := G_i^* G_j, \quad i, j = 1, \ldots, b,$$

where $A_{ij}$ is, generally, a rectangular $n_i \times n_j$ matrix. Hence, all matrices $A^{(k)} = [G^{(k)}]^* G^{(k)}$ have the same induced partition as the initial matrix $A$.

Typically, in the $k$-th block-step of the algorithm, the block-pivoting strategy chooses a pair of indices $(i, j)$, such that $1 \leq i < j \leq b$. This pair determines two block-columns, or the so-called *block-pivot* submatrix in the current working matrix $G^{(k-1)}$, that will be transformed in this block-step. Like before, we shall omit the superscript $(k-1)$ when referring to the individual blocks of $G^{(k-1)}$ and $A^{(k-1)}$. The block-pivot submatrix of $G^{(k-1)}$ is then defined by

$$G_P^{(k-1)} := [G_i, G_j], \qquad (3.5)$$

and it has $n^{(k-1)} = n_i + n_j$ individual columns.

$$A_P^{(k-1)} := \begin{bmatrix} A_{ii} & A_{ij} \\ A_{ij}^* & A_{jj} \end{bmatrix} = [G_P^{(k-1)}]^* G_P^{(k-1)} = \begin{bmatrix} G_i^* G_i & G_i^* G_j \\ G_j^* G_i & G_j^* G_j \end{bmatrix}. \qquad (3.6)$$

In both algorithms (full block and block-oriented), in addition to the two block-column pivot submatrix $G_P^{(k-1)}$ from (3.5), the single block-column pivot $G_P^{(k-1)} := [G_i]$ is employed. In this case, the corresponding square block-pivot submatrix of $A^{(k-1)}$ is

$$A_P^{(k-1)} := [G_P^{(k-1)}]^* G_P^{(k-1)} = G_i^* G_i.$$

Block-partitions of $G^{(k-1)}$ and $A^{(k-1)}$ are illustrated in Figure 1. Both block-pivot submatrices are shown in color.

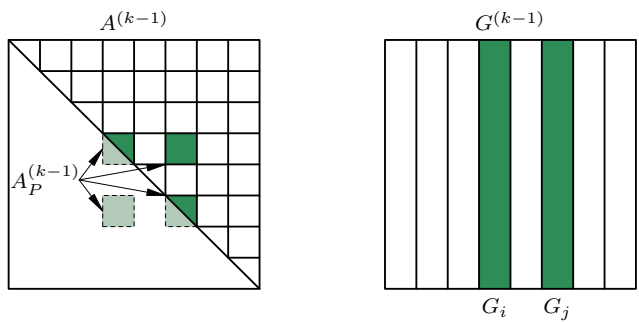

Figure 1: Block-columns $[G_i, G_j]$ in $G^{(k-1)}$ generate a square block $A_P^{(k-1)}$ in $A^{(k-1)} = [G^{(k-1)}]^* G^{(k-1)}$.

Note that $G_P^{(k-1)}$ is always of full column rank, so $A_P^{(k-1)}$ is certainly positive definite. Also, $G_P^{(k-1)}$ is one of many possible "Cholesky-like" factors of $A_P^{(k-1)}$. To speed up the process, instead of "tall" factor $G_P^{(k-1)}$, the Cholesky factor $R_P^{(k-1)}$ of $A_P^{(k-1)}$ will be used. Transformations applied to $R_P^{(k-1)}$ will be accumulated and applied to $G_P^{(k-1)}$ in BLAS 3 fashion by block matrix multiplication.

The only difference between the full block and block-oriented algorithms lies in the transformation of block-pivot submatrices $A_P^{(k-1)}$ from (3.6).

### 3.3 Full block algorithm

The main idea behind the full block algorithm is to work with block-columns in a similar fashion as the non-blocked one-sided algorithm works with individual columns.

In the non-blocked algorithm, the pivot strategy chooses a pair of columns that generates a matrix pair $(A^{(k-1)}, J^{(k-1)})$ of order 2 in (3.6), which will be diagonalized in step $k$. This simply means that the off-diagonal element $A_{ij}$ will be annihilated in that step.

In the full block algorithm, the so-called outer, or block-pivot strategy chooses a pair of block-columns for transformation in this step. Likewise, they generate a matrix pair $(A^{(k-1)}, J^{(k-1)})$ in (3.6), but the order of this pair is, generally, greater than 2.

The non-blocked algorithm can now be mimicked in two different ways.
1. We can annihilate only the off-diagonal block $A_{ij}$, i.e., block-diagonalize the pivot submatrix $A^{(k-1)}$. In the end, this strategy produces a block-diagonal matrix. The final step consists of the diagonalization of these diagonal blocks. Unfortunately, it can be shown (see [18]) that this approach is only linearly convergent to the block-diagonal form.
2. Therefore, in each step of the algorithm, the whole pivot submatrix $A^{(k-1)}$ is diagonalized. Even though this strategy seems less efficient at first sight, it is quadratically convergent. If the block sizes are chosen to be sufficiently small, so that most of the required data resides in cache, the matrix $A^{(k-1)}$ can be efficiently diagonalized by the non-blocked one-sided $J$-Jacobi algorithm.

In this approach, after the first full block-sweep, the diagonal blocks of the working matrix $A$ remain diagonal throughout the whole process. Thus, it pays off to do a preprocessing step, before the block-sweep loop, in which all diagonal blocks are diagonalized. Afterwards, in each step we have to diagonalize a pivot submatrix of the following form

$$A^{(k-1)} = \begin{bmatrix} \Lambda_{ii} & A_{ij} \\ A_{ij}^* & \Lambda_{jj} \end{bmatrix}, \tag{3.7}$$

where $\Lambda_{ii}$ and $\Lambda_{jj}$ are diagonal matrices.

The one-sided $J$-Jacobi algorithm works on a factor of $A^{(k-1)}$. This process will be the fastest possible, if it works on a square factor, not on a "tall" one. A suitable square factor $R^{(k-1)}$ can be computed from (3.6) in two ways:
1. by the QR factorization of $[G_i, G_j]$, which is a "tall" factor of $A^{(k-1)}$, or
2. by the Cholesky factorization of $A^{(k-1)}$.

Due to the special structure (3.7) of $A^{(k-1)}$, the second choice turns out to be much faster, as it requires only one matrix multiplication to form the off-diagonal block $A_{ij}$.

Moreover, the structure of $A^{(k-1)}$ in (3.7) can also be exploited to efficiently compute the Cholesky factorization

$$A^{(k-1)} = [R^{(k-1)}]^* R^{(k-1)}, \qquad R^{(k-1)} = \begin{bmatrix} R_{ii} & R_{ij} \\ & R_{jj} \end{bmatrix}. \tag{3.8}$$

Since the diagonal blocks of $A^{(k-1)}$ in (3.7) are already diagonal, it is easy to see that

$$R_{ii} = \sqrt{\Lambda_{ii}} \tag{3.9}$$

$$R_{ij} = (R_{ii})^{-1} A_{ij} \tag{3.10}$$

$$R_{jj}^* R_{jj} = \Lambda_{jj} - R_{ij}^* R_{ij}. \tag{3.11}$$

The final block $R_{jj}$ can be computed by the Cholesky factorization of the right-hand side in (3.11). An alternative is obtained by writing

$$\Lambda_{jj} - R_{ij}^* R_{ij} = \begin{bmatrix} R_{ij}^* & \sqrt{\Lambda_{jj}} \end{bmatrix} \operatorname{diag}(-I_{ii}, I_{jj}) \begin{bmatrix} R_{ij} \\ \sqrt{\Lambda_{jj}} \end{bmatrix} := \widehat{G} \widehat{J} \widehat{G}^*,$$

and we can use the hyperbolic QR factorization of $\widehat{G}^*$ according to $\widehat{J}$ to compute $R_{jj}$ (see [28] for details). Both factorizations are accurate [29], but once again, the Cholesky factorization is faster.

The diagonalization of $A^{(k-1)}$ is now equivalent to the orthogonalization of columns of the Cholesky factor $R^{(k-1)}$ from (3.8). This is done by using the "in-cache" non-blocked one-sided hyperbolic Jacobi method to compute the HSVD of the upper triangular matrix $R^{(k-1)}$, and yields

$$R^{(k-1)} W^{(k)} = U^{(k)} \Sigma^{(k)}, \tag{3.12}$$

or

$$R^{(k-1)} = U^{(k)} \Sigma^{(k)} [V^{(k)}]^*, \quad V^{(k)} = [W^{(k)}]^{-*}.$$

From (3.8) and (3.12) we obtain

$$[W^{(k)}]^* A^{(k-1)} W^{(k)} = [W^{(k)}]^* [R^{(k-1)}]^* R^{(k-1)} W^{(k)} = [\Sigma^{(k)}]^2,$$

so, $W^{(k)}$ is the matrix which simultaneously diagonalizes both $A^{(k-1)}$ and $J^{(k-1)}$. Note that $W^{(k)}$ can be computed in several different ways (see [18]), but numerical tests show that the explicit accumulation of all used $J$-rotations is the best option, both regarding the speed, and the accuracy. A careful analysis shows that during the first few sweeps of the block method, the accumulation of transformations (on a micro-level) is relatively slow, but in the last few sweeps, when the number of rotations applied becomes progressively smaller, the accumulation becomes faster and faster.

The final and most time-consuming part of the block algorithm is the so-called update (2.9) of block-columns in $G^{(k-1)}$. Post-multiplication of $[G_i, G_j]$ by $W^{(k)}$ is done by the BLAS 3 routine xGEMM. However, since xGEMM cannot overwrite the original factor matrix $[G_i, G_j]$, we use an additional $n \times n_0$ array as a workspace, where $n_0 = 2 \max_{i=1,\ldots,p} n_i$. To avoid unnecessary copying of the pivot blocks, updated columns are not moved back to their original positions in the matrix. We just keep track of their current position, and in a post-processing step, all block-columns are re-permuted back to their original positions.

Since two- and three-level algorithms are based on both sequential non-blocked and blocked algorithms, we will present them here for the sake of completeness. In the first auxilliary subroutine, `Jacobi_Cycle` (see Algorithm 3.1) either off-diagonal elements of $A$, or the off-diagonal block $A_{ij}$ of $A$, are annihilated. The routine is a slightly modified version of the inner loop (for Hermitian matrices instead of symmetric) of Algorithm 1 from [32].

---

**Algorithm 3.1:** Single sweep of the hyperbolic Jacobi algorithm

`Jacobi_Cycle(G, J, D, W, m, n_i, n_j, diag_bl);`
**Description**: If $diag\_bl = true$ then $G$ consists of a single block column, i.e.,
$G = G_i \in \mathbb{C}^{m \times n_i}$ and $J = J_{ii}$. Else, $G$ consists of two block columns,
$G = [G_i, G_j] \in \mathbb{C}^{m \times (n_i + n_j)}$ and $J = \text{diag}(J_{ii}, J_{jj})$. The diagonal of $A = G^*G$
is stored in vector $D$, i.e., $D = \text{diag}(A)$. Transformations applied to $G$ are
accumulated in matrix $W$.

**begin**
  **if** $diag\_bl$ **then**    // parameters for triangular part of $A$
  | $start\_s = 2$;  $final\_s = n_i$;
  **else**    // parameters for rectangular block $A_{ij} = G_i^* G_j$
  | $start\_s = n_i + 1$;  $final\_s = n_i + n_j$;
  **end if**
  **for** $s = start\_s$ **to** $final\_s$ **do**
    **if** $diag\_bl$ **then** $final\_r = s - 1$ **else** $final\_r = n_i + n_j$;
    **for** $r = 1$ **to** $final\_r$ **do**
      $a_{rs} = g_r^* g_s$;  $hyp = -J_{rr} \cdot J_{ss}$;
        // compute parameters for modified rotation defined by (2.3)
      $\eta = \text{sign}(\text{re}\, a_{rs}) \cdot |a_{rs}|$;  $e^{i\alpha} = a_{rs}/\eta$;   // if $A$ is real, $\eta = a_{rs}$
      $\vartheta = hyp \cdot (a_{ss} - hyp \cdot a_{rr})/(2 \cdot \eta)$;  $t = \text{sign}(\vartheta)/(|\vartheta| + \sqrt{\vartheta^2 + hyp})$;
      $h = \sqrt{1 + hyp \cdot t^2}$;  $cs = 1/h$;  $sn = cs \cdot t$;
        // update diagonal of $D$
      $d_{rr} = d_{rr} + hyp \cdot t \cdot \eta$;  $d_{ss} = d_{ss} + t \cdot \eta$;
        // update columns $r$ and $s$ of $G$ and $W$
      $f = e^{i\alpha} \cdot g_r$;  $g_r = cs \cdot f - hyp \cdot sn \cdot g_s$;  $g_r = sn \cdot f + cs \cdot g_s$;
      $f = e^{i\alpha} \cdot w_r$;  $w_r = cs \cdot f - hyp \cdot sn \cdot w_s$;  $w_r = sn \cdot f + cs \cdot w_s$;
    **end for**
  **end for**
**end**

---

If the `Jacobi_Cycle` algorithm is repeatedly applied to the off-diagonal elements of $A$ until convergence, the obtained algorithm is the one-sided hyperbolic algorithm (called `Jacobi_Diagonalization` in Algorithm 3.2). To speed up the calculation (see [32] for details) the routine also computes the diagonal of $A$. Finally, in Algorithm 3.3 we present the sequential full block algorithm.

---

**Algorithm 3.2:** One-sided hyperbolic Jacobi diagonalization algorithm

`Jacobi_Diagonalization(G, J, D, W, m, n, max_it)`;
**Description**: Diagonalization of pair $(A, J)$, where $A = G^*G$, $G \in \mathbb{C}^{m \times n}$ by the one-sided hyperbolic Jacobi algorithm, until either convergence, or maximal number of iterations $max\_it$ are reached. Vector $D$ keeps the diagonal of $A$. Transformations applied to $G$ are accumulated in matrix $W$.

**begin**
    $iter = 0$;   $W = I_n$;
    **repeat**
        $iter = iter + 1$;
        **for** $i = 1$ **to** $n$ **do** $d_{ii} = g_i^* g_i$;   // in each cycle initialize the diagonal of $D$
        `Jacobi_Cycle(G, J, D, W, m, n, 0, true)`;
    **until** *(convergence) or (iter $\geq$ max_it)*;
**end**

---

**Algorithm 3.3:** Sequential full block algorithm

`Full_Block(G, J, D, W, m, n, max_it)`;
**Description**: Diagonalization of pair $(A = G^*G, J)$, $G \in \mathbb{C}^{m \times n}$ by the full block Jacobi algorithm.

**begin**
    partition $G$ into almost equally-sized column blocks $G_i$, $i = 1, \ldots, b$;
    corresponding partition of $J$ is $J = \text{diag}(J_{11}, \ldots, J_{bb})$;
    $W = I_n$;   $iter = 0$;
    **repeat**
        $iter = iter + 1$;
        **for** $i = 1$ **to** $b$ **do**   // diagonal block pre-processing
            $A_P = G_i^* G_i$;
            Cholesky factorization of $A_P = R_{ii}^* R_{ii}$;
            `Jacobi_Diagonalization(R_ii, J_ii, D_ii, W_ii, n_i, n_i, max_it)`;
            $G_i = G_i \cdot W_{ii}$;   $W_i = W_i \cdot W_{ii}$
        **end for**
        **for** $j = 2$ **to** $b$ **do**   // diagonalization of pivot block $A_P$
            **for** $i = 1$ **to** $j - 1$ **do**
                $A_P = [G_i\ G_j]^*[G_i\ G_j]$;   $J_P = \text{diag}(J_{ii}, J_{jj})$;   $n_P = n_i + n_j$;
                Cholesky factorization of $A_P = R^* R$;
                `Jacobi_Diagonalization(R, J_P, D_P, W_P, n_P, n_P, max_it)`;
                $[G_i\ G_j] = [G_i\ G_j] \cdot W_P$;   $[W_i\ W_j] = [W_i\ W_j] \cdot W_P$
            **end for**
        **end for**
    **until** *(convergence) or (iter $\geq$ max_it)*;
**end**

---

## 3.4 Block-oriented algorithm

The block-oriented algorithm is a simple rearrangement of each sweep of the one-sided $J$-Jacobi algorithm in a cache-aware manner.

Typically, in each sweep, the off-diagonal elements of diagonal blocks $A_{ii}$ in (3.6) are annihilated first. Then, for each pivot submatrix $A^{(k-1)}$, all elements in the off-diagonal block $A_{ij}$ are annihilated once, thus completing a sweep.

As in the full block case, the algorithm works on a factor of $A^{(k-1)}$, and again it is advantageous to have a square factor $R^{(k-1)}$. Since $A^{(k-1)}$ has no special structure in the block-oriented algorithm, the factor $R^{(k-1)}$ is computed by using the standard Cholesky factorization.

The accumulation of $W^{(k)}$ and the update of block-columns in $G^{(k-1)}$ is done in exactly the same way as in the full block algorithm.

---
**Algorithm 3.4:** Sequential block oriented-algorithm

---
`Block_Oriented`($G$, $J$, $D$, $W$, $m$, $n$, $max\_it$);
**Description**: Diagonalization of pair $(A = G^*G, J)$, $G \in \mathbb{C}^{m \times n}$ by the block-oriented Jacobi algorithm.
**Assumption**: If this function is used as a part of a parallel three-level block-oriented algorithm, partition of $G$ should respect column boundaries of $G_{i\_blk}$ and $G_{j\_blk}$.

**begin**
    partition $G$ into almost equally-sized column blocks $G_i$, $i = 1, \ldots, b$;
    corresponding partition of $J$ is $J = \text{diag}(J_{11}, \ldots, J_{bb})$;
    $W = I_n$;    $iter = 0$;
    **repeat**
        $iter = iter + 1$;
        **for** $i = 1$ **to** $b$ **do**    // single pass through diagonal blocks
            $W_P = I_{n_i}$;    $A_P = G_i^* G_i$;
            Cholesky factorization of $A_P = R_{ii}^* R_{ii}$;
            `Jacobi_Diagonalization`($R_{ii}$, $J_{ii}$, $D_{ii}$, $W_{ii}$, $n_i$, $n_i$, 1);
            $G_i = G_i \cdot W_{ii}$;    $W_i = W_i \cdot W_{ii}$
        **end for**
        **for** $j = 2$ **to** $b$ **do**
            **for** $i = 1$ **to** $j - 1$ **do**    // single pass through off-diagonal blocks
                $W_P = I_{n_i + n_j}$;    $J_P = \text{diag}(J_{ii}, J_{jj})$;    $A_P = [G_i \, G_j]^* [G_i \, G_j]$;
                Cholesky factorization of $A_P = R^* R$;
                `Jacobi_Cycle`($R$, $J_P$, $D_P$, $W_P$, $m$, $n_i$, $n_j$, false);
                $[G_i \, G_j] = [G_i \, G_j] \cdot W_P$;    $[W_i \, W_j] = [W_i \, W_j] \cdot W_P$
            **end for**
        **end for**
    **until** *(convergence)* or *(iter $\geq$ max_it)*;
**end**

---

## 3.5 Sequential testing

The first tests were made on a single core processor, (Pentium 660J, 3.6GHz, running Windows 64-bit), in double precision. We were also using Intel Fortran compiler and the BLAS and LAPACK routines from the Intel Math Kernel Library.

For block algorithms, a provably convergent strategy (see [17, 18]): the block column-cyclic at the block level, and the ordinary column-cyclic inside each block is used. The same strategy is used in parallel three-level algorithms inside each process.

The eigenvalues of symmetric test matrices are randomly generated in the given range, by using the LAPACK routine `xLARND`. The matrices are then generated by the routine `xLAGSY` for generation of symmetric matrices, based on random reflectors. Test matrices have various orders, from 500 to 4000 in steps of 500, and the block-columns $G_i$ consist of 8–128 individual columns.

Our tests show that the non-blocked Jacobi algorithm runs faster than its block counterparts for matrices of order less than 1000. For small orders, a big portion of the matrix resides in the cache memory, while book-keeping by the block algorithms takes too much time. For matrices of order $1000 \leq n < 2500$, the block-oriented algorithm is the fastest one. Finally, for matrices of order $n \geq 2500$, the fastest algorithm is the full block algorithm.

The speedups of the full block and the block-oriented algorithms, compared to the non-blocked algorithm, for $n = 4000$ are displayed in Figure 2.

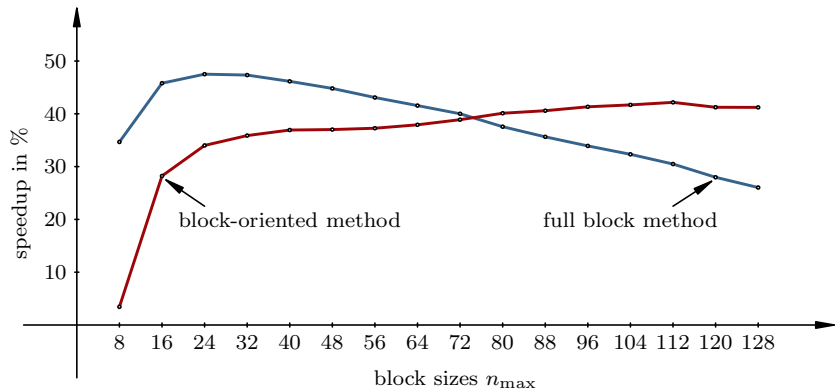

Figure 2: Speedup of sequential blocked algorithms for matrix size $n = 4000$, compared to non-blocked algorithm, single core Pentium 660J processor.

For matrices of order $n$ between 1500 and 3500, the speedups are somewhat smaller, but have the same shape as in Figure 2.

The second test was made on a single core of Intel Xeon E5420 processor (single core of test machine $A$, see later). We were quite surprised that the block-oriented algorithm is faster than the full block algorithm for all tested matrix sizes (see Figure 3).

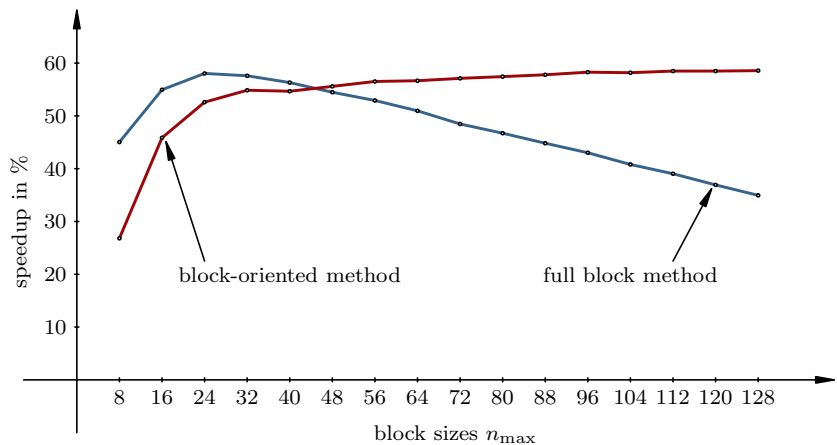

Figure 3: Speedup of sequential blocked algorithms for matrix size $n = 4000$, compared to non-blocked algorithm, one core of Xeon processor.

As neither of the sequential blocking strategies is clearly better than the other, we shall continue to use both of them in parallel algorithms, as well.

## 4 Parallel two-level algorithms

The idea of blocking can also be used to parallelize the Jacobi algorithm, but with an entirely different goal, to ensure data independency. This will be done simply by changing the block-pivot strategy. Except for the communication between processes, everything else is inherited from sequential block algorithms, so the computational procedure is essentially the same.

The basic idea of parallelization is to partition $G$ and $J$ as in (3.1), and map this partition to processes. From now on, $p$ will denote the number of available processes. Then, the reasonable number of blocks in partition (3.1) is $2p$. One can use $2p - 1$ blocks, as shown in Figure 4, but in this case, one process has less work to do than the others.

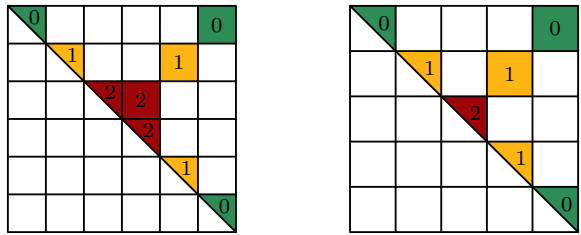

Figure 4: Independent blocks for $p = 3$, for even ($n = 6$) and odd ($n = 5$) number of blocks. If $n = 5$, the workload of process 2 is approximately half of the workload of processes 0 or 1.

As we have already said, there are two standard parallel pivoting strategies that can be used as block-pivoting strategies in parallel Jacobi algorithms.
1. The modulus pivot strategy is known to be convergent, but it uses one parallel step more per sweep than is minimally required.
2. The round-robin strategy is appealing because it uses the minimal number of parallel steps per sweep. Unfortunately, nothing is known about its convergence properties.

We decided to implement both strategies and compare them in practice.

The ring or one-dimensional torus is a natural topology of the parallel processes in our algorithm. To minimize the communication in each parallel step, only one block-column of the current working matrix $G$ is sent to one neighboring process, and the other block-column is received from the other neighboring process. This simple and uniform communication pattern is in fact a simultaneous one-step cyclic shift in the ring of processes. The network bandwidth is, therefore, fully utilized. The amount of transmitted data remains (almost) constant in the whole run, thus making the communication overhead predictable and scalable, if the network speed increases.

Due to blocking, memory access in the communication routines is contiguous, and the relative impact of the network latency effectively diminishes as the problem size increases.

Our first goal is to map the block-columns $[G_i, G_j]$ and the corresponding diagonal blocks of $J$ to individual processes in such a way to ensure data independency between processes and a natural communication pattern between them.

The mapping to processes is inherited from the modulus strategy, with a slight shift with respect to the original modulus strategy presented in [21], so that our modulus strategy begins on the main antidiagonal of the matrix. This is the initial layout for both parallel pivoting strategies. More precisely, the initial order of block-columns of $G$ with respect to processes is $(1, 2p), (2, 2p-1), \ldots, (p, p+1)$, i.e., block-columns $(q+1, 2p-q)$ reside in process $q$, for $q = 0, \ldots, p - 1$.

In the MPI terminology, $q$ is usually called a *rank* of a process, and $p$ is the total number of processes (denoted by *nproc* in all code segments). Each process, by knowing only its own rank and the number of blocks $nbl := 2 \cdot nproc$, independently computes indices $i\_blk$, and $j\_blk$ of its first pivot pair in a sweep.

## 4.1 Modulus pivoting strategy

Each block-sweep consists of a certain number of parallel steps according to the chosen block-pivot strategy. A single parallel step is a collection of all non-overlapping transformations executed in parallel (in all processes). Figure 5 shows the block layouts for the modulus pivot strategy in odd and even sweeps.

In each sweep of the modulus strategy, the blocks denoted by ▫ are visited twice, and we have two possibilities:

**Algorithm 4.1:** Initialization of modulus and modified round-robin strategy

Initialize (*rank*, *nbl*, *ip*, *jp*, *i_blk*, *j_blk*);
**Description**: At the start of a sweep, the algorithm for process *rank* computes indices of the initial column blocks (*i_blk*, *j_blk*). The auxilliary pair (*ip*, *jp*) is used for determination of the pivot indices in the next steps.

**begin**
$\quad$ $ip = rank + 1;\quad i\_blk = ip;\quad jp = nbl - rank;\quad j\_blk = jp;$
**end**

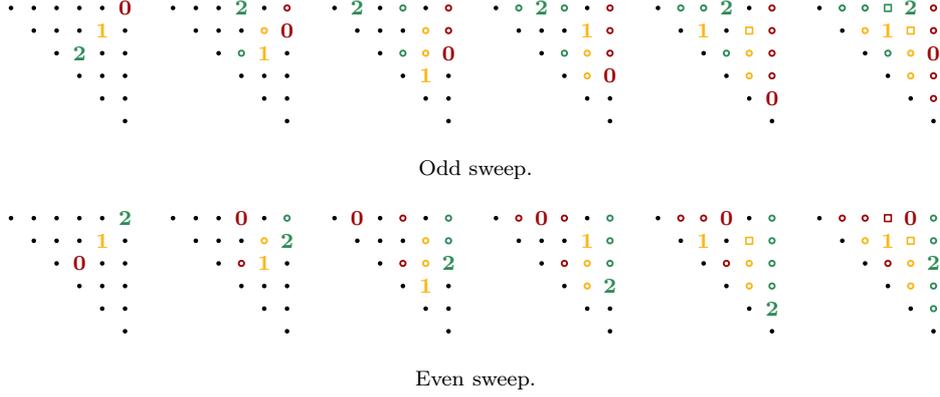

Figure 5: Modulus strategy for $p = 3$ processes and $2p = 6$ blocks. A number denotes the index of a process where the block (two block-columns) resides, • denotes still unvisited blocks, ○ denotes already processed blocks, while □ denotes the blocks that will be visited twice.

(a) transform the block twice, as in the so-called quasi-cyclic strategies,
(b) just do nothing, and wait for the next pair of blocks.

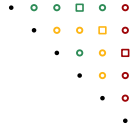

Figure 6: For $n = 3$ and $p = 6$, the blocks visited twice are denoted by □.

In our implementation, these blocks are transformed twice. Note that doubly transformed blocks are always at the same place, occupying the $n$-th superdiagonal. Figure 6 shows the blocks that are transformed twice in the first sweep.

Subroutine `MMStep` in Algorithm 4.2 calculates the next pivot indices for process *rank* according to the modified modulus strategy. By using the sweep-number parity, *nsweep*, it determines the next process-sender *rcv_rnk*, which sends its block-column *rcv_blk* to *rank*, and the next process-receiver *snd_rnk*, which receives block-column *snd_blk* from *rank*.

## 4.2 Round-robin pivoting strategy

The standard round-robin tournament strategy has to be modified to ensure that each process sends and receives only one block-column in each parallel step. The strategy is best described by grouping the indices of the block-columns in two horizontal groups, where the indices above

**Algorithm 4.2:** A step of a modified modulus strategy

MMStep (*rank*, *nbl*, *nproc*, *nsweep*, *ip*, *jp*, *i_blk*, *j_blk*, *snd_rnk*, *snd_blk*, *rcv_rnk*, *rcv_blk*);
**Description**: If a step is not the first one in a sweep, process with rank *rank* determines the indices of the next pivot block $(i\_blk, j\_blk)$.

**begin**
    **if** $(ip + jp) > nbl$ **then**
        $snd\_blk = i\_blk$;    $ip = ip + 1$;
        **if** $ip = jp$ **then**
            $ip = ip - nbl/2$;    $jp = ip$;
        **end if**
        $i\_blk = ip$;    $rcv\_blk = i\_blk$;
    **else**
        $snd\_blk = j\_blk$;    $jp = jp + 1$;    $j\_blk = jp$;    $rcv\_blk = j\_blk$;
    **end if**
    **if** $(nsweep \bmod 2) > 0$ **then**
        $snd\_rnk = (nproc + rank - 1) \bmod nproc$;    $rcv\_rnk = (nproc + rank + 1) \bmod nproc$;
    **else**
        $snd\_rnk = (nproc + rank + 1) \bmod nproc$;    $rcv\_rnk = (nproc + rank - 1) \bmod nproc$;
    **end if**
**end**

each other belong to the columns in the same process.

In the standard strategy, the indices are rotated counter-clockwise by one to obtain the next layout. In the modified strategy, the indices are rotated in the clockwise direction by $\lfloor (2p-1)/2 \rfloor$. This process is illustrated in Figure 7.

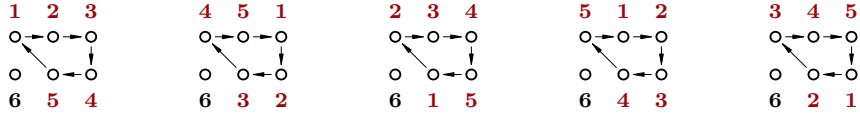

Figure 7: Modified round-robin strategy for $p = 3$.

It is worth to mention that our ordering is "preferable" to Brent and Luk's ordering [5] since it is not working on elements reduced in the previous step (see [14, page 455] for details), at least for small number of processes $p$. Unfortunately, the proof whether or not is modified round-robin strategy convergent seems to be uneasy.

The block layouts in odd and even sweeps for the modified round-robin strategy are shown in Figure 8.

Subroutine RLStep in Algorithm 4.3 calculates the next indices of blocks $(i\_blk, j\_blk)$ which will reside in process *rank*, according to the modified round-robin strategy. All used variables have the same meaning as in subroutine MMStep. The only difference is variable *swflag*. The value *swflag* = 1 means that the blocks in a process *rank* are placed in reverted positions, i.e., $i\_blk > j\_blk$, and should be renumerated such that $i\_blk < j\_blk$.

It should be said that, despite its optimality with respect to the number of parallel steps, in most of our tests, the round-robin turned out to be slower than the modulus strategy.

## 4.3 Two-level algorithms

In a given parallel step, each process has is own local pair of block-columns $G^{[q]} := [G_i, G_j]$, where $q$ is the index of the process, that generates a "local" block-pivot matrix. To simplify the notation, this block-pivot matrix will be denoted by $A^{[q]}$, and appropriate $J$ by $J^{[q]}$. So,

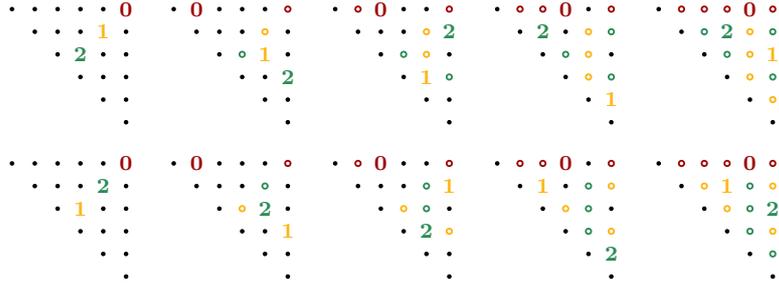

Figure 8: Modified round-robin strategy for $p = 3$ processes and $2p = 6$ blocks. A number denotes the index of a process where the block (two block-columns) resides, • denotes still unvisited blocks, and ○ denotes already processed blocks.

---

**Algorithm 4.3:** A step of the modified round-robin strategy

RLStep (*rank*, *nbl*, *nproc*, *nsweep*, *ip*, *jp*, *i_blk*, *j_blk*, *snd_rnk*, *snd_blk*, *rcv_rnk*, *rcv_blk*, *swflag*);

**Description**: If the step is not the first one in the sweep, process with rank *rank* determines the indices of the next pivot block (*i_blk*, *j_blk*).

**begin**
  *swflag* = 0;
  **if** $(ip + jp) > nbl$ **then**
    **if** $jp < nbl$ **then**
      *snd_blk* = *i_blk*;   *ip* = *ip* + 1;
      **if** $ip < jp$ **then**
        *i_blk* = *ip*;   *rcv_blk* = *i_blk*;
      **else**
        *ip* = *jp*;   *i_blk* = *j_blk*;   *swflag* = 1;   *jp* = *nbl*;   *j_blk* = *jp*;
        *rcv_blk* = *j_blk*;
      **end if**
    **else**
      **if** $ip > (nbl/2)$ **then**
        *snd_blk* = *i_blk*;   *ip* = *ip* − *nbl*/2 + 1;   *i_blk* = *ip*;   *rcv_blk* = *i_blk*;
      **else**
        *snd_blk* = *j_blk*;   *jp* = *ip* + 1;   *j_blk* = *jp*;   *rcv_blk* = *j_blk*;
      **end if**
    **end if**
  **else**
    *snd_blk* = *j_blk*;   *jp* = *jp* + 1;   *j_blk* = *jp*;   *rcv_blk* = *j_blk*;
  **end if**
  **if** (*nsweep* mod 2) > 0 **then**
    *snd_rnk* = (*nproc* + *rank* − 1) mod *nproc*;   *rcv_rnk* = (*nproc* + *rank* + 1) mod *nproc*;
  **else**
    *snd_rnk* = (*nproc* + *rank* + 1) mod *nproc*;   *rcv_rnk* = (*nproc* + *rank* − 1) mod *nproc*;
  **end if**
**end**

---

as in (3.6), we have

$$A^{[q]} = \begin{bmatrix} A_{ii} & A_{ij} \\ A_{ij}^* & A_{jj} \end{bmatrix} = [G^{[q]}]^* G^{[q]} = [G_i, G_j]^* [G_i, G_j] = \begin{bmatrix} G_i^* G_i & G_i^* G_j \\ G_j^* G_i & G_j^* G_j \end{bmatrix}. \quad (4.1)$$

In two-level algorithms, the first step is different than the others.

1. The matrix pair $(A^{[q]}, J^{[q]})$, for each process $q$ is diagonalized by the $J^{[q]}$-unitary trans-

formations, but the special formulae (3.9)–(3.11) cannot be used for the Cholesky factorization, since the diagonal blocks of $A^{[q]}$ are not yet diagonal.

2. In the block-oriented case the full sweep of $J^{[q]}$-unitary transformations is applied to the off-diagonal elements of $A^{[q]}$.

The other steps in a sweep consist of the following computational tasks in each process $q$:

1. local matrix pair $(A^{[q]}, J^{[q]})$ is
    (a) diagonalized in the full block case, but note that now $A^{[q]}$ has diagonal blocks already diagonalized, and (3.9)–(3.11) can be used,
    (b) in the block-oriented case, only elements of block $A_{ij}$ are annihilated.

    In both cases the local transformation matrix $W^{[q]}$ is accumulated throughout this process, and

2. a local pair of block-columns is updated by $W^{[q]}$, i.e., $[G^{[q]}]' := [G'_i, G'_j] = G^{[q]} W^{[q]}$.

This is followed by a communication step to interchange block-columns between processes.

By using previous two algorithms, and the communication, we can easily write both two-level algorithms (see Algorithm 5.2).

# 5 Parallel three-level algorithms

If we compare the optimal block sizes (which depend on the order $n$ of a matrix) for the sequential full block and block-oriented algorithms in Figures 2 and 3, we see that the optimal sizes are somewhere between 8 and 48 for the full block algorithm, and between 24 and 128 for the block-oriented algorithm.

This suggests that, as soon as the block size in each process becomes sufficiently large, it is advisable to use a sequential blocked algorithm for local diagonalization, instead of its non-blocked counterpart. This motivates us to construct the "doubly" blocked or three-level versions of parallel $J$-Jacobi algorithms and compare their efficiency to the non-blocked (two-level) parallel algorithms.

## 5.1 Three-level full block algorithm

The idea behind the three-level full block algorithm is to use the sequential blocked algorithm with optimal sequential block size inside each process. Note that inside each process we have two block-columns, each with $n$ rows and approximately $n_b$ columns, with $n_b \ll n$.

The first step of the algorithm is to "shorten" these two "tall" local columns $[G_i, G_j]$. These columns are multiplied to obtain the local block-pivot matrix $A^{[q]}$ as in (4.1), and the Cholesky factorization of $A^{[q]}$ gives a significantly smaller factor $R^{[q]}$ of $A^{[q]}$—the order of $R^{[q]}$ is approximately $2n_b$. The sequential blocked one-sided Jacobi algorithm is then used to compute the HSVD of $R^{[q]}$ with respect to $J^{[q]}$.

In the final step, the accumulated transformation matrix $W^{[q]}$ that diagonalizes the local pair $(A^{[q]}, J^{[q]})$ is applied to update the block-columns $[G_i, G_j]$. Note that if the second level (in cache, or inner) method is chosen to be the full block method, it is not necessary to take care whether the inner block partition respects the boundary between two "outer" blocks $[G_i, G_j]$, or not.

## 5.2 Three-level block-oriented algorithm

The basic idea is almost the same as in the three-level full block algorithm, but we have to take care that the inner block partition respects the boundary between the two "outer" local blocks $[G_i, G_j]$.

Here, the distinction has to be made between two different types of sub-blocks in the local block-pivot matrix $A^{[q]}$ from (4.1).

1. The diagonal sub-blocks $A_{ii}$ pose no problem, as they are diagonalized only once in each block-sweep, and the corresponding block-column $G_i$ can be locally partitioned in the usual manner.
2. To deal with, generally, rectangular off-diagonal block $A_{ij} = G_i^* G_j$, each of the local blocks $G_i$ and $G_j$ have to be partitioned individually, as in (3.1). This correctly defines an "inner" partition of the block $R_{ij}$ in the Cholesky factorization of $A^{[q]}$, and this block is used by the sequential block-oriented one-sided $J$-Jacobi algorithm. This partitioning procedure is given in more detail in Algorithm 5.1.

---

**Algorithm 5.1:** Off-diagonal step of parallel three-level block-oriented algorithm

`Off_Diagonal(G, J, D, W, m, n_r, n_s)` **Description**: Single annihilation pass in two-level blocked manner of the off-diagonal block of $A_{rs} = G_r^* G_s$ of $A = G^* G$, where $G = [G_r\ G_s]$ by the block-oriented Jacobi algorithm.

**Assumption**: Block $G_i$ has $n_i$ columns.

**begin**
  partition $G_r$ into almost equally-sized column blocks $G_i$, $i = 1, \ldots, b_r$;
  partition $G_s$ into almost equally-sized column blocks $G_i$, $i = b_r + 1, \ldots, b$;
  $b = b_r + b_s$;   $n = n_r + n_s$;   $W = I_n$;
  corresponding partition of $J$ is $J = \mathrm{diag}(J_{11}, \ldots, J_{bb})$;
  **for** $j = b_r + 1$ **to** $b$ **do**    // `single pass through` $A_{rs}$
    **for** $i = 1$ **to** $b_r$ **do**
      $W_P = I_{n_i+n_j}$;   $A_P = [G_i\ G_j]^* [G_i\ G_j]$;   $J_P = \mathrm{diag}(J_{ii}, J_{jj})$;
      Cholesky factorization of $A_P = R^* R$;
      `Jacobi_Cycle(`$R$, $J_P$, $D_P$, $W_P$, $m$, $n_i$, $n_j$, *false*`)`;
      $[G_i\ G_j] = [G_i\ G_j] \cdot W_P$;   $[W_i\ W_j] = [W_i\ W_j] \cdot W_P$
    **end for**
  **end for**
**end**

---

## 5.3 The core of the parallel algorithms

The main similarities and differences of the two-level and three-level parallel algorithms are summarized at the top level by the following description.

The factor $G$ is distributed evenly to $p$ processes as $2p$ block-columns of $n$ rows, and approximately $\frac{n}{2p}$ columns each, where the number of columns varies by at most one between all block-columns. A pair $(q+1, 2p-q)$ of block-columns is assigned to a process $q$, for $0 \leq q < p$. Processes are organized into a one-dimensional torus, where each process has exactly two neighbors to communicate with. In our implementation, each process is bound to a single processor core.

As a first step in each cycle, the *full block* (F) strategy diagonalizes diagonal blocks $A^{[q]}$ in parallel to obtain diagonal blocks in diagonal form (3.7) for the rest of the sweep. Note that this could be done only in the first sweep, but we prefer occasional re-diagonalization for the sake of the accuracy of the special Cholesky factorization (3.9)–(3.11). In the *block-oriented* (B) algorithms, in the first step of each sweep the off-diagonal elements of the whole $A^{[q]}$ are annihilated, while in the rest of the sweep, only elements of the off-diagonal block $A_{ij}$ in (4.1) are annihilated.

The algorithm works in sweeps, until convergence is detected. Convergence test is twofold: either no rotations in the previous sweep were applied, or rotation angles were all below the predefined quadratic convergence threshold.

Each sweep makes $2p$ (modulus) or $2p-1$ (round-robin) steps. In each step $p$ independent pairs of block-columns $[G_i, G_j]$ are transformed in parallel. The choice of pairs is directed by the parallel pivot strategy. Each pair of block-columns is acted upon at least once in a sweep.

The anatomy of a step for the process $q$ follows:
1. A "small" matrix $A^{[q]} = [G_i, G_j]^*[G_i, G_j]$ is formed from the two block-columns, $G_i$ and $G_j$, local to process $q$.
2. The Cholesky factorization is performed to obtain $A^{[q]} = [R^{[q]}]^* R^{[q]}$, with $R^{[q]}$ square. In the *full block* strategy, $R^{[q]}$ is computed by using (3.9)–(3.11).
3. The columns of $R^{[q]}$ are transformed by the hyperbolic one-sided Jacobi algorithm:
   (a) In the *full block* strategy, the hyperbolic singular-value decomposition of $R^{[q]}$ is computed, such that $R^{[q]} W^{[q]} = U^{[q]} \Sigma^{[q]}$, i.e. columns of $R^{[q]}$ are orthogonalized.
   (b) In the *block-oriented* (B) strategy, the off-diagonal elements of $A^{[q]}$ are annihilated exactly once by using congruences. Instead of applying transformations to both sides of $A^{[q]}$, i.e., $[A^{[q]}]' = [W^{[q]}]^* A^{[q]} W^{[q]}$, they are applied to the right-hand side of $R^{[q]}$, i.e., $[R^{[q]}]' = R^{[q]} W^{[q]}$.
   
   All the transformations applied are accumulated in the matrix $W^{[q]}$ in both cases. If $R^{[q]}$ is large enough, the orthogonalization (annihilation) is performed by the blocked, otherwise by the non-blocked one-sided hyperbolic Jacobi algorithm. This is the only difference between three-level (3B, 3F) and two-level (2B, 2F) algorithms.
4. The accumulated transformations $W^{[q]}$ are applied to the right-hand side of $[G_i, G_j]$, i.e., $[G'_i, G'_j] = [G_i, G_j] W^{[q]}$.
5. According to the parallel pivot strategy, one block-column is chosen to be sent to the one of the neighboring processes, and is simultaneously replaced by another block-column received from the other neighbor.

This description, written in pseudo-code as Algorithm 5.2, shows how simple the efficient parallel Jacobi algorithm can be built from the essential routines: the non-blocked and blocked sequential algorithms together with simple communication routines.

# 6 Numerical testing

Our parallel Jacobi algorithms are implemented in Fortran, using Intel MKL for sequential BLAS and LAPACK routines, HDF5 libraries for storing data, and synchronous point-to-point and collective communication routines of MPI-1 from Open MPI distribution. There are no inherent constraints in the code on hardware platform, cluster configuration or matrix dimensions, save physical and method-induced ones. Test processes are 64-bit and single-threaded. Each process is bound to its own processor core and communicates with others through MPI stack only.

## 6.1 The test machines

Test machine $A$ is a five-blade cluster. Each blade has two Intel Xeon E5420 processors, 16 GB DDR2 RAM, and a gigabit Ethernet network card connected to a switch. The processors run at 2.5 GHz, with 1333 MHz FSB, $2 \times 6$ MB level-2 cache memory, and SSE 4.1 instruction set available. Each processor has four cores with independent $32 + 32$ kB level-1 caches. The cores of a processor share level-2 cache in pairs. Cache memory at both levels runs at the core speed.

For the initial testing two other parallel machines, $B$ and $C$, were used. Machine $B$ is a single SMP server with NUMA architecture, Sun Fire X4600 M2, with 8 dual-core AMD Opteron 8220 CPUs (core speed 2.8 GHz, 1 MB level-2 cache per core), and 48 GB DDR2/667 ECC RAM. Machine $C$ is a cluster of 16 blades, each blade consisting of two dual-core Intel Xeon 5150 CPUs (core speed 2.66 GHz, $32 + 32$ kB level-1 cache per core, and 4 MB level-2 cache shared between the cores), 8 GB DDR2 RAM, and a gigabit Ethernet network card connected to a switch.

**Algorithm 5.2:** Parallel block Jacobi algorithms

Parallel_Jacobi($G$, $J$, $m$, $n$);

**Description**: Diagonalization of a Hermitian matrix $H$, given in factored form $H = GJG^*$, $G \in \mathbb{C}^{m \times n}$ via the one-sided hyperbolic SVD of $G$.

**Assumption**: $G$ is divided into $nbl$ block-columns, and $J$ into $nbl$ diagonal blocks. In each process the first block-column is denoted by index $i\_blk$, and the second by $j\_blk$.

**begin**

  Initialize ($rank$, $nbl$, $ip$, $jp$, $i\_blk$, $j\_blk$);

  receive/read appropriate $G^{[q]} = [G_{i\_blk} \ G_{j\_blk}]$, $J^{[q]} = \text{diag}(J_{i\_blk,i\_blk}, J_{j\_blk,j\_blk})$;

  **repeat**

    $first\_step = true$;    $n^{[q]} = n_{i\_blk} + n_{j\_blk}$;

    // compute square factor of matrix $A^{[q]}$;

    $A^{[q]} = [G^{[q]}]^* [G^{[q]}]$;

    Cholesky factorization of $A^{[q]} = [R^{[q]}]^* R^{[q]}$;

    // call apropriate sequential algorithm for each block;

    **switch** *algorithm* **of**

      **case** *2F*:

        Jacobi_Diagonalization($R^{[q]}$, $J^{[q]}$, $D^{[q]}$, $W^{[q]}$, $n^{[q]}$, $n^{[q]}$, $max\_it$);

      **endcase**

      **case** *2B*:

        **if** *first_step* **then**

          Jacobi_Diagonalization($R^{[q]}$, $J^{[q]}$, $D^{[q]}$, $W^{[q]}$, $n^{[q]}$, $n^{[q]}$, 1);

          $first\_step = false$;

        **else**

          Jacobi_Cycle($R^{[q]}$, $J^{[q]}$, $D^{[q]}$, $W^{[q]}$, $n^{[q]}$, $n_{i\_blk}$, $n_{j\_blk}$, $false$);

        **end if**

      **endcase**

      **case** *3F*:

        Full_Block($R^{[q]}$, $J^{[q]}$, $D^{[q]}$, $W^{[q]}$, $n^{[q]}$, $n^{[q]}$, $max\_it$)

      **endcase**

      **case** *3B*:

        **if** *first_step* **then**

          Block_Oriented($R^{[q]}$, $J^{[q]}$, $D^{[q]}$, $W^{[q]}$, $n^{[q]}$, $n^{[q]}$, 1);    $first\_step = false$;

        **else**

          Off_Diagonal($R^{[q]}$, $J^{[q]}$, $D^{[q]}$, $W^{[q]}$, $n^{[q]}$, $n_{i\_blk}$, $n_{j\_blk}$);

        **end if**

      **endcase**

    **endswitch**

    $G^{[q]} = G^{[q]} \cdot W^{[q]}$;

    // calculate send/receive indices, for example by MMStep;

    MMStep ($rank$, $nbl$, $nproc$, $nsweep$, $ip$, $jp$, $i\_blk$, $j\_blk$, $snd\_rnk$, $snd\_blk$, $rcv\_rnk$, $rcv\_blk$);

    send $G_{snd\_blk}$, $D_{snd\_blk}$, and $J_{snd\_blk,snd\_blk}$ to $snd\_rnk$;

    receive $G_{rcv\_blk}$, $D_{rcv\_blk}$, and $J_{rcv\_blk,rcv\_blk}$ from $rcv\_rnk$;

    exchange convergence informations;

  **until** *convergence*;

  eigenvectors: normalized columns of $G^{[q]}$;

  eigenvalues: squared norms of columns of $G^{[q]}$ multiplied by signs from $J^{[q]}$;

**end**

All machines are running 64-bit GNU/Linux (Red Hat on $A$, Debian on $B$ and $C$).

Due to operating system's process scheduler tendency to occasionally reschedule a process to a different core or even processor socket from one it was running on, what induces unwanted cache spills and possibly issues of non-local memory accesses NUMA architectures are suffering from, a simple scheduling protocol is devised to bind processes permanently during a test run to processor cores, in a one-to-one manner. Very similar features have recently been implemented in Open MPI version 1.3, and since then have been used in our testing instead.

Initial tests on machine $B$ have exhibited severe timing inconsistencies between test runs of the same algorithm over identical data sets, differing by a factor of up to 2.5. Repeated identical runs of processes pinned to cores our protocol (or, lately, Open MPI) dedicated for them, however, differ by no more than 1–2% on machine $B$, and by just a few seconds on machines $A$ and $C$.

## 6.2 The test results

A production setup of the three-level $J$-Jacobi algorithms should begin with the careful profiling of the corresponding (full block or block-oriented) sequential blocked $J$-Jacobi algorithm, to reveal the optimal inner block size for the particular hardware. Our testing has shown that there exists a global optimum over all inner block sizes for the given machine and the problem size. The optimal block size is expected to be a slowly growing function of the problem size on any modern processor architecture.

Similarly, the decision between the full block and block-oriented three level algorithms is platform and problem type dependant. The testing suggests that for smaller problems and in an environment where interprocess communication is relatively slow (e.g., Ethernet) the full block algorithm is a method of choice. For larger problems and with a fast communication available (e.g., local shared memory) the block-oriented algorithm wins. It is up to the implementor, however, to correctly identify the adequate algorithm for the targeted machines and spectra of problem sizes.

Our test results show that the modified round-robin strategy is always slower than the modulus strategy. On machine $A$, it is about 3% slower on average, while on machine $B$, for very large matrices, it can be even two times slower. Therefore, we shall present more detailed test results only for the modulus strategy.

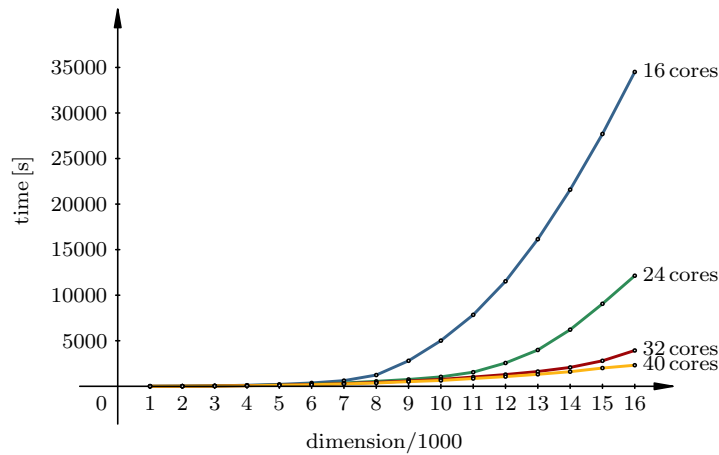

Figure 9: Modulus strategy: measured time for two-level full block method (2F), for $p = 16, \ldots, 40$ cores.

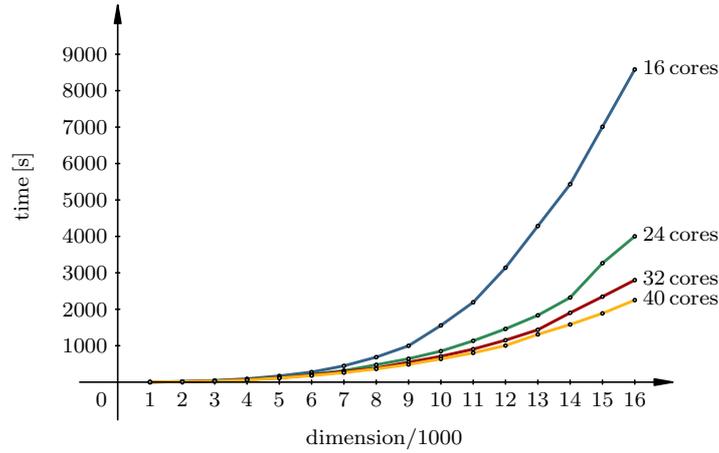

Figure 10: Modulus strategy: measured time for two-level block-oriented method (2B), for $p = 16, \ldots, 40$ cores.

As we expected, a significant speedup is obtained if we compare three-level versus two-level algorithms. In the case of modulus pivot strategy, the ratio between two- and three-level algorithm is approximately 30% for sufficiently large matrices. From Figure 13 it is clear that the usage of three-level method is advisable (on machine $A$) if the order of matrix is 9000 or more.

In the case of block-oriented strategy (see Figure 14), the ratio between two- and three-level methods (on machine $A$) is approximately 55% for large matrices.

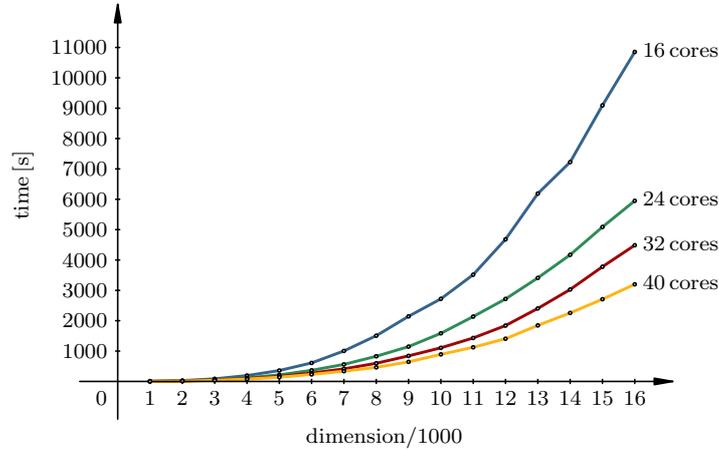

Figure 11: Modulus strategy: measured time for three-level full block method (3F), for $p = 16, \ldots, 40$ cores.

The previous two graphs show us only that by usage of blocking we obtained the significant speedup, but they do not show the true relation between speeds of these two methods. Finally, if we compare times for various algorithms (see Figures 9–12), we can conclude that the modulus version of the three-level block-oriented algorithm is the fastest one.

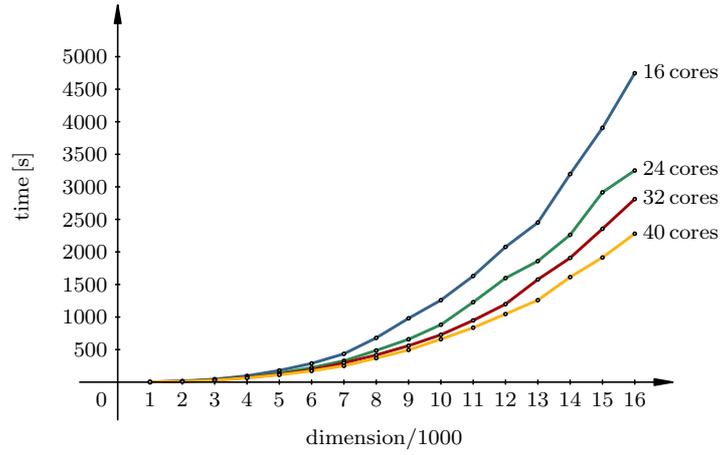

Figure 12: Modulus strategy: measured time for three-level block-oriented method (3B), for $p = 16, \ldots, 40$ cores.

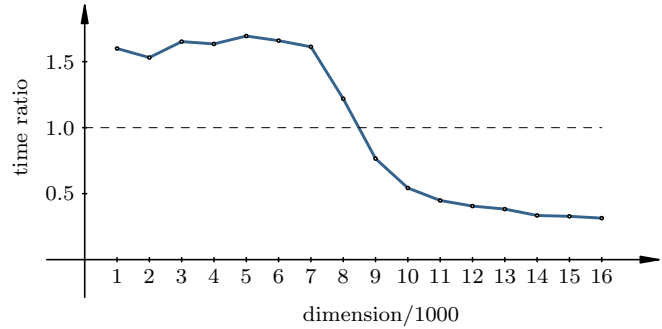

Figure 13: The ratio of times for three-level full block and two-level full block method (modulus pivot strategy), $p = 16$ cores.

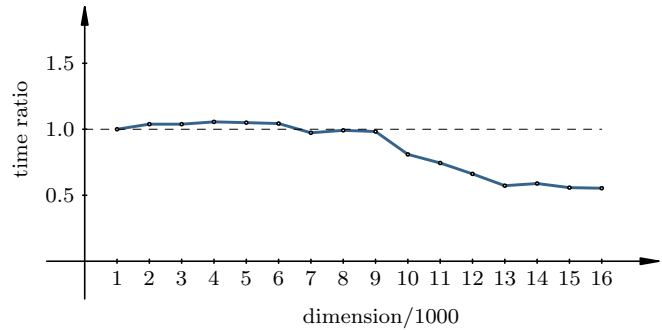

Figure 14: The ratio of times for three-level block-oriented and two-level block-oriented method (modulus pivot strategy), $p = 16$ cores.

## 6.3 A simple time model

From the theoretical viewpoint, complexity of the sequential Jacobi algorithm is proportional to $c_n \cdot n^3$. The algorithm is almost ideally paralellizable, and it can be divided into $p$ almost equally sized processes. Therefore the expected execution time should be

$$\text{Time}(n,p) = c\frac{n^3}{p},$$

where $c$ is a parameter that depends on the chosen algorithm, $n$, $p$, and the speed of the interprocess communication.

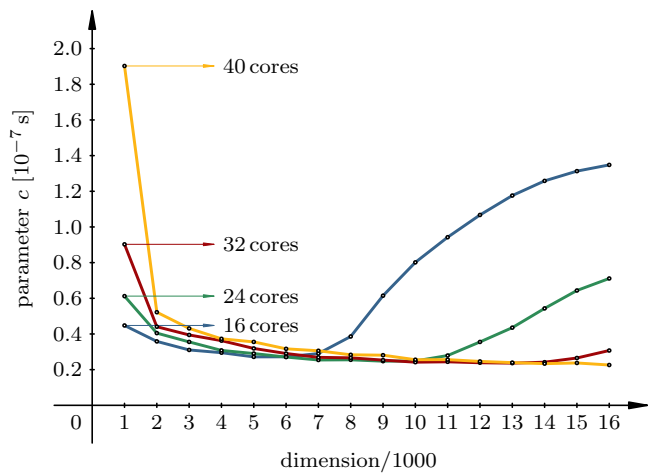

Figure 15: Modulus strategy: parameter $c$ for 2F, $\text{Time}(n,p) = c\frac{n^3}{p}$.

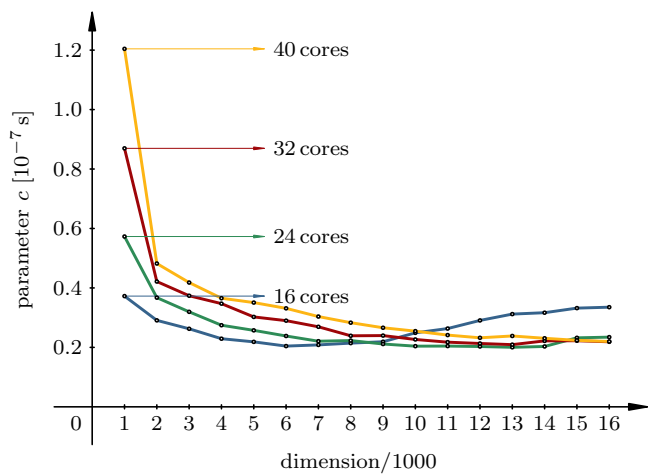

Figure 16: Modulus strategy: parameter $c$ for 2B, $\text{Time}(n,p) = c\frac{n^3}{p}$.

We really cannot expect that $c$ is a constant. For example, the number of sweeps is a slowly growing function of $n$. Nevertheless, the experiments show that $c$ is almost constant in the middle range, when neither communication, nor the block size inside each core are the dominant factors. Figures 15–18 show that phenomenon.

Finally, for a more accurate time prediction model we need significantly larger test space and a precise model that captures subtle properties of the processor, cache, memory and network behavior.

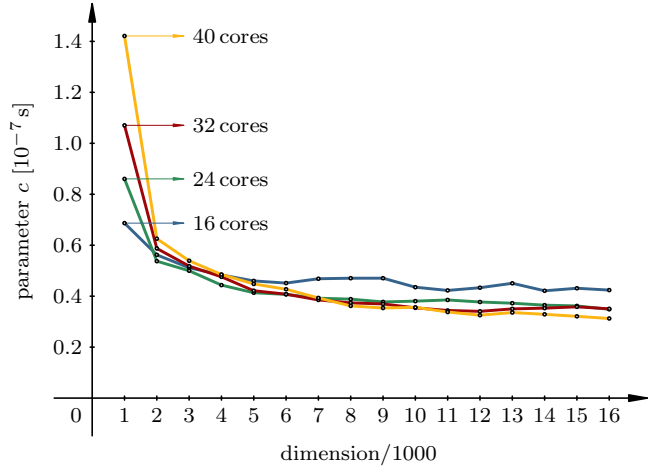

Figure 17: Modulus strategy: parameter $c$ for 3F, $\text{Time}(n,p) = c\frac{n^3}{p}$.

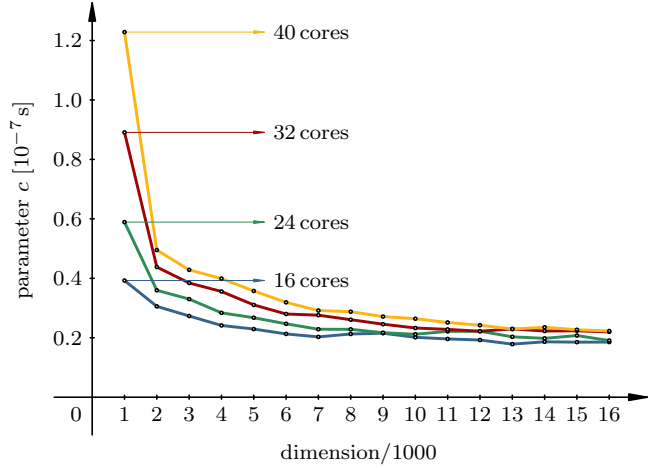

Figure 18: Modulus strategy: parameter $c$ for 3B, $\text{Time}(n,p) = c\frac{n^3}{p}$.

# Conclusion

In this paper we have introduced the two- and three-level parallel hyperbolic Jacobi algorithms. These types of organizations of parallel sweeps can make diagonalization process significantly faster for the large enough matrices. The convergence of the algorithms (with modulus block strategy and cyclic by rows/columns inside blocks) is a direct consequence of the convergence of the row/column cyclic block algorithms [17, 18].

Our bunch of algorithms permits fine tuning for intricacies of the targeted computer architecture. The choice of block-oriented or full block algorithm depends on the communication speed between processes—the full block algorithms are suitable for the slower speed, and block-oriented for the faster. Once the crossings between two- and three-level algorithms are found on the particular architecture, the three-level algorithms can be additionally finely tuned—the size of inner blocks should be determined for optimal performance.

The testing results show almost perfect scalability (see Figures 15–18), in a sense that $c$ is really a constant for the middle ranges of $n$ and $p$.

## Acknowledgements

We would like to express special thanks to our friend and colleague prof. Vjeran Hari who gave us an idea to describe the full block algorithm in the spirit of the *J*-Jacobi algorithm.

We would also like to thank prof. Dorian Marjanović, prof. Nenad Bojčetić, and Tomislav Rašeta, from the Faculty of Mechanical Engineering and Naval Architecture, University of Zagreb, for providing a parallel machine for the initial tests of our algorithms. The testing continued at GFZ German Research Centre for Geosciences, Helmholtz-Centre Potsdam, due to generosity of prof. Mioara Mandea and dr. Vincent Lesur. We are also deeply indebted to Vedran Dunjko and Robert Bakarić from Ruđer Bošković Institute, Zagreb, who allowed us to complete the testing on their cluster.